\documentstyle[12pt,amssymb,epic,eepic,fullpage]{amsart}
\def\draft{n}
\theoremstyle{plain}

\newtheorem{proposition}{Proposition}[section]
\newtheorem{lemma}[proposition]{Lemma}

\theoremstyle{definition}

\newtheorem{lemmadefinition}[proposition]{Lemma-Definition}

\theoremstyle{remark}
\newtheorem{example}[proposition]{Example}

\newtheorem{remark}[proposition]{Remark}

\def\printname#1{
	\if\draft y
		\smash{\makebox[0pt]{\hspace{-0.5in}
			\raisebox{8pt}{\tt\tiny #1}}}
	\fi
}

\newcommand{\mathmode}[1]{$#1$}

\newlength{\standardunitlength}
\catcode`\@=11
\long\def\@makecaption#1#2{%
    \vskip 10pt
    \setbox\@tempboxa\hbox{
      \small\sf{\bfcaptionfont #1. }\ignorespaces #2}%
    \ifdim \wd\@tempboxa >\captionwidth {%
        \rightskip=\@captionmargin\leftskip=\@captionmargin
        \unhbox\@tempboxa\par}%
      \else
        \hbox to\hsize{\hfil\box\@tempboxa\hfil}%
    \fi}
\font\bfcaptionfont=cmssbx10 scaled \magstephalf
\newdimen\@captionmargin\@captionmargin=2\parindent
\newdimen\captionwidth\captionwidth=\hsize
\catcode`\@=12

\newcommand{\sign}{\operatorname{sign}}

\def\lbl#1{\label{#1}\printname{#1}}

\theoremstyle{plain}
\newtheorem{statement}{Statement}

\newcommand{\tr}{\operatorname{tr}}
\newcommand{\Wst}{W_{\widetilde{so(3)}}}
\newcommand{\EEPIC}[2]{
	\setlength{\unitlength}{#2\standardunitlength}
        \begin{array}{c}  \hspace{-1.7mm}
                \raisebox{-8pt}{#1}
                \hspace{-1.9mm}
        \end{array}
}

\def\BWConvention
{{
\begin{picture}(6250,1164)(0,-10)
\thicklines
\path(309.114,937.496)(237.000,837.000)(347.923,891.737)
\put(430.269,609.115){\arc{597.609}{4.0090}{7.7476}}
\path(1062,537)(1362,537)
\path(1242.000,507.000)(1362.000,537.000)(1242.000,567.000)
\path(4662,537)(4962,537)
\path(4842.000,507.000)(4962.000,537.000)(4842.000,567.000)
\path(462,612)(462,12)
\path(12,1062)(462,612)(912,1062)
\path(2038,612)(2038,12)
\path(2188,612)(2188,12)
\path(2113,687)(2563,1137)
\path(2188,612)(2638,1062)
\path(1663,1137)(2113,687)
\path(1588,1062)(2038,612)
\path(3909.114,937.496)(3837.000,837.000)(3947.923,891.737)
\put(4030.269,609.115){\arc{597.609}{4.0090}{7.7476}}
\path(4063,612)(4063,12)
\put(3762,462){\makebox(0,0)[lb]{\smash{{\mathmode{\scriptstyle -}}}}}
\path(3613,1062)(4063,612)(4513,1062)
\path(6163,1137)(5938,912)(5938,837)
\path(6238,1062)(6013,837)(5938,837)
\path(5263,1137)(5488,912)(5488,837)
	(5413,837)(5188,1062)
\path(5638,12)(5638,312)(5713,387)
	(5788,312)(5788,12)
\path(5713,387)(5788,462)(5788,612)
	(5938,762)(5938,837)(5863,837)
	(5713,687)(5563,837)(5488,837)
	(5488,762)(5638,612)(5638,462)
	(5713,387)(5713,387)
\put(2038,87){\makebox(0,0)[lb]{\smash{{\mathmode{\scriptstyle w}}}}}
\put(5638,87){\makebox(0,0)[lb]{\smash{{\mathmode{\scriptstyle w}}}}}
\put(5653,537){\makebox(0,0)[lb]{\smash{{\mathmode{\scriptstyle b}}}}}
\put(162,462){\makebox(0,0)[lb]{\smash{{\mathmode{\scriptstyle +}}}}}
\end{picture}
}}
\def\DoubleEdges
{{
\begin{picture}(3646,645)(0,-10)
\thicklines
\path(2412,411)(3612,411)
\path(2412,261)(3612,261)
\path(1662,336)(1962,336)
\path(1842.000,306.000)(1962.000,336.000)(1842.000,366.000)
\put(87,411){\makebox(0,0)[lb]{\smash{{\mathmode{\scriptstyle a}}}}}
\path(12,336)(1212,336)
\put(1062,411){\makebox(0,0)[lb]{\smash{{\mathmode{\scriptstyle b}}}}}
\put(2487,36){\makebox(0,0)[lb]{\smash{{\mathmode{\scriptstyle i}}}}}
\put(3462,486){\makebox(0,0)[lb]{\smash{{\mathmode{\scriptstyle k}}}}}
\put(3462,36){\makebox(0,0)[lb]{\smash{{\mathmode{\scriptstyle l}}}}}
\put(2487,486){\makebox(0,0)[lb]{\smash{{\mathmode{\scriptstyle j}}}}}
\end{picture}
}}
\def\Flip
{{
\begin{picture}(8473,2537)(0,-10)
\thicklines
\put(1833,1614){\ellipse{1800}{1800}}
\put(1458,1914){\ellipse{150}{150}}
\put(2208,1914){\ellipse{150}{150}}
\path(333,1614)(933,1614)
\path(2733,1614)(3333,1614)
\path(1833,1914)(1833,1314)
\put(6634.000,1727.500){\arc{975.000}{3.5364}{5.8884}}
\put(6634,1615){\ellipse{1800}{1800}}
\put(6259,1315){\ellipse{150}{150}}
\put(1833.000,1501.500){\arc{975.000}{0.3948}{2.7468}}
\put(7009,1315){\ellipse{150}{150}}
\put(5583,39){\makebox(0,0)[lb]{\smash{{{\rm\scriptsize rest of the graph}}}}}
\path(5134,1615)(5734,1615)
\path(7534,1615)(8134,1615)
\path(6634,1315)(6634,1915)
\path(3933,1314)(4533,1314)
\path(4413.000,1284.000)(4533.000,1314.000)(4413.000,1344.000)
\path(333,1614)	(286.682,1606.427)
	(246.880,1598.778)
	(184.189,1582.376)
	(139.653,1563.036)
	(108.000,1539.000)

\path(108,1539)	(77.311,1491.352)
	(55.867,1429.144)
	(41.740,1365.614)
	(33.000,1314.000)

\path(33,1314)	(27.048,1271.794)
	(21.889,1223.971)
	(17.524,1171.332)
	(13.952,1114.681)
	(11.174,1054.820)
	(9.190,992.551)
	(7.999,928.677)
	(7.601,864.000)
	(7.998,799.323)
	(9.188,735.449)
	(11.172,673.180)
	(13.950,613.319)
	(17.521,556.668)
	(21.887,504.029)
	(27.046,456.206)
	(33.000,414.000)

\path(33,414)	(41.745,362.386)
	(55.871,298.856)
	(77.312,236.648)
	(108.000,189.000)

\path(108,189)	(155.648,158.311)
	(217.856,136.868)
	(281.386,122.740)
	(333.000,114.000)

\path(333,114)	(377.580,109.236)
	(437.753,107.648)
	(475.883,108.045)
	(520.549,109.236)
	(572.628,111.221)
	(633.000,114.000)

\path(3340,1614)	(3386.318,1606.427)
	(3426.120,1598.778)
	(3488.811,1582.376)
	(3533.347,1563.036)
	(3565.000,1539.000)

\path(3565,1539)	(3595.689,1491.352)
	(3617.132,1429.144)
	(3631.260,1365.614)
	(3640.000,1314.000)

\path(3640,1314)	(3645.952,1271.794)
	(3651.111,1223.971)
	(3655.476,1171.332)
	(3659.048,1114.681)
	(3661.826,1054.820)
	(3663.810,992.551)
	(3665.001,928.677)
	(3665.399,864.000)
	(3665.002,799.323)
	(3663.812,735.449)
	(3661.828,673.180)
	(3659.050,613.319)
	(3655.479,556.668)
	(3651.113,504.029)
	(3645.954,456.206)
	(3640.000,414.000)

\path(3640,414)	(3631.255,362.386)
	(3617.129,298.856)
	(3595.688,236.648)
	(3565.000,189.000)

\path(3565,189)	(3517.352,158.311)
	(3455.144,136.868)
	(3391.614,122.740)
	(3340.000,114.000)

\path(3340,114)	(3295.420,109.236)
	(3235.247,107.648)
	(3197.117,108.045)
	(3152.451,109.236)
	(3100.372,111.221)
	(3040.000,114.000)

\path(5133,1614)	(5086.682,1606.427)
	(5046.880,1598.778)
	(4984.189,1582.376)
	(4939.653,1563.036)
	(4908.000,1539.000)

\path(4908,1539)	(4877.311,1491.352)
	(4855.868,1429.144)
	(4841.740,1365.614)
	(4833.000,1314.000)

\path(4833,1314)	(4827.048,1271.794)
	(4821.889,1223.971)
	(4817.524,1171.332)
	(4813.952,1114.681)
	(4811.174,1054.820)
	(4809.190,992.551)
	(4807.999,928.677)
	(4807.601,864.000)
	(4807.998,799.323)
	(4809.188,735.449)
	(4811.172,673.180)
	(4813.950,613.319)
	(4817.521,556.668)
	(4821.887,504.029)
	(4827.046,456.206)
	(4833.000,414.000)

\path(4833,414)	(4841.745,362.386)
	(4855.871,298.856)
	(4877.312,236.648)
	(4908.000,189.000)

\path(4908,189)	(4955.648,158.311)
	(5017.856,136.868)
	(5081.386,122.740)
	(5133.000,114.000)

\path(5133,114)	(5177.580,109.236)
	(5237.753,107.648)
	(5275.883,108.045)
	(5320.549,109.236)
	(5372.628,111.221)
	(5433.000,114.000)

\path(8140,1614)	(8186.318,1606.427)
	(8226.120,1598.778)
	(8288.811,1582.376)
	(8333.347,1563.036)
	(8365.000,1539.000)

\path(8365,1539)	(8395.689,1491.352)
	(8417.133,1429.144)
	(8431.260,1365.614)
	(8440.000,1314.000)

\path(8440,1314)	(8445.952,1271.794)
	(8451.111,1223.971)
	(8455.476,1171.332)
	(8459.048,1114.681)
	(8461.826,1054.820)
	(8463.810,992.551)
	(8465.001,928.677)
	(8465.399,864.000)
	(8465.002,799.323)
	(8463.812,735.449)
	(8461.828,673.180)
	(8459.050,613.319)
	(8455.479,556.668)
	(8451.113,504.029)
	(8445.954,456.206)
	(8440.000,414.000)

\path(8440,414)	(8431.255,362.386)
	(8417.129,298.856)
	(8395.688,236.648)
	(8365.000,189.000)

\path(8365,189)	(8317.352,158.311)
	(8255.144,136.868)
	(8191.614,122.740)
	(8140.000,114.000)

\path(8140,114)	(8095.420,109.236)
	(8035.247,107.648)
	(7997.117,108.045)
	(7952.451,109.236)
	(7900.372,111.221)
	(7840.000,114.000)

\put(783,39){\makebox(0,0)[lb]{\smash{{{\rm\scriptsize rest of the graph}}}}}
\end{picture}
}}
\def\ThickG
{{
\begin{picture}(2509,2231)(0,-10)
\thicklines
\put(1137.500,970.500){\arc{1629.801}{1.4940}{3.9921}}
\put(1412.500,970.500){\arc{1629.801}{5.4327}{7.9308}}
\put(1275,983){\ellipse{2100}{1950}}
\path(1350,158)(1350,608)(1200,758)
	(1200,983)(1050,1133)(1050,1283)(675,1658)
\path(600,1583)(975,1208)(1125,1208)(1275,1058)
\path(1950,1583)(1575,1208)(1425,1208)(1275,1058)
\path(1200,158)(1200,608)(1350,758)
	(1350,983)(1500,1133)(1500,1283)(1875,1658)
\put(1200,2108){\makebox(0,0)[lb]{\smash{{\mathmode{\scriptstyle u}}}}}
\put(1650,1133){\makebox(0,0)[lb]{\smash{{\mathmode{\scriptstyle i}}}}}
\put(1275.000,970.500){\arc{1825.000}{3.9948}{5.4299}}
\put(1425,1433){\makebox(0,0)[lb]{\smash{{\mathmode{\scriptstyle j}}}}}
\put(525,683){\makebox(0,0)[lb]{\smash{{\mathmode{\scriptstyle p}}}}}
\put(975,1433){\makebox(0,0)[lb]{\smash{{\mathmode{\scriptstyle k}}}}}
\put(675,1133){\makebox(0,0)[lb]{\smash{{\mathmode{\scriptstyle l}}}}}
\put(1425,533){\makebox(0,0)[lb]{\smash{{\mathmode{\scriptstyle n}}}}}
\put(0,608){\makebox(0,0)[lb]{\smash{{\mathmode{\scriptstyle q}}}}}
\put(1950,683){\makebox(0,0)[lb]{\smash{{\mathmode{\scriptstyle r}}}}}
\put(2325,608){\makebox(0,0)[lb]{\smash{{\mathmode{\scriptstyle s}}}}}
\put(1200,1658){\makebox(0,0)[lb]{\smash{{\mathmode{\scriptstyle t}}}}}
\put(900,533){\makebox(0,0)[lb]{\smash{{\mathmode{\scriptstyle m}}}}}
\end{picture}
}}
\def\ThickVertex
{{
\begin{picture}(2937,1221)(0,-10)
\thicklines
\path(675,612)(675,12)
\path(600,687)(1050,1137)
\path(675,612)(1125,1062)
\path(150,1137)(600,687)
\path(75,1062)(525,612)
\path(2325,12)(2325,312)(2475,462)(2475,612)
\path(2475,12)(2475,312)(2325,462)(2325,612)
\path(2850,1137)(2625,912)(2625,762)(2475,612)
\path(2925,1062)(2700,837)(2550,837)(2400,687)
\path(1949,1136)(2174,911)(2174,761)(2324,611)
\path(1874,1061)(2099,836)(2249,836)(2399,686)
\path(525,612)(525,12)
\put(2776,762){\makebox(0,0)[lb]{\smash{{\mathmode{\scriptstyle i}}}}}
\put(1425,537){\makebox(0,0)[lb]{\smash{{\mathmode{-}}}}}
\put(2551,1062){\makebox(0,0)[lb]{\smash{{\mathmode{\scriptstyle j}}}}}
\put(2101,1062){\makebox(0,0)[lb]{\smash{{\mathmode{\scriptstyle k}}}}}
\put(1801,762){\makebox(0,0)[lb]{\smash{{\mathmode{\scriptstyle l}}}}}
\put(2551,162){\makebox(0,0)[lb]{\smash{{\mathmode{\scriptstyle n}}}}}
\put(975,762){\makebox(0,0)[lb]{\smash{{\mathmode{\scriptstyle i}}}}}
\put(750,1062){\makebox(0,0)[lb]{\smash{{\mathmode{\scriptstyle j}}}}}
\put(300,1062){\makebox(0,0)[lb]{\smash{{\mathmode{\scriptstyle k}}}}}
\put(0,762){\makebox(0,0)[lb]{\smash{{\mathmode{\scriptstyle l}}}}}
\put(750,162){\makebox(0,0)[lb]{\smash{{\mathmode{\scriptstyle n}}}}}
\put(225,162){\makebox(0,0)[lb]{\smash{{\mathmode{\scriptstyle m}}}}}
\put(2025,162){\makebox(0,0)[lb]{\smash{{\mathmode{\scriptstyle m}}}}}
\end{picture}
}}
\def\TwoFamilies
{{
\begin{picture}(6624,2646)(0,-10)
\thicklines
\put(5256.026,1352.132){\arc{2718.434}{5.6899}{7.6862}}
\put(5525.316,1058.904){\arc{2716.345}{3.9456}{5.3959}}
\put(5290.684,1058.904){\arc{2716.344}{4.0289}{5.4792}}
\put(1208,1212){\ellipse{2400}{2400}}
\path(308,2037)(1208,1212)(2108,2037)
\path(1208,1212)(1208,12)
\path(3008,1212)(3608,1212)
\path(3488.000,1182.000)(3608.000,1212.000)(3488.000,1242.000)
\path(4583,2037)(5333,1287)
\path(5483,1287)(6233,2037)
\path(5333,1287)(5483,12)
\path(5483,1287)(5333,12)
\put(1433,1662){\makebox(0,0)[lb]{\smash{{\mathmode{\scriptstyle 1}}}}}
\put(5559.974,1352.132){\arc{2718.433}{1.7386}{3.7348}}
\put(608,1437){\makebox(0,0)[lb]{\smash{{\mathmode{\scriptstyle 2}}}}}
\put(5108,237){\makebox(0,0)[lb]{\smash{{\mathmode{\scriptstyle 1}}}}}
\put(1283,537){\makebox(0,0)[lb]{\smash{{\mathmode{\scriptstyle 3}}}}}
\put(233,687){\makebox(0,0)[lb]{\smash{{\mathmode{\scriptstyle 1}}}}}
\put(2033,687){\makebox(0,0)[lb]{\smash{{\mathmode{\scriptstyle 2}}}}}
\put(1208,2112){\makebox(0,0)[lb]{\smash{{\mathmode{\scriptstyle 3}}}}}
\put(5633,1662){\makebox(0,0)[lb]{\smash{{\mathmode{\scriptstyle 1}}}}}
\put(4808,1437){\makebox(0,0)[lb]{\smash{{\mathmode{\scriptstyle 2}}}}}
\put(6233,687){\makebox(0,0)[lb]{\smash{{\mathmode{\scriptstyle 2}}}}}
\put(4508,687){\makebox(0,0)[lb]{\smash{{\mathmode{\scriptstyle 1}}}}}
\put(5933,2487){\makebox(0,0)[lb]{\smash{{\mathmode{\scriptstyle 2}}}}}
\put(5708,2112){\makebox(0,0)[lb]{\smash{{\mathmode{\scriptstyle 1}}}}}
\put(5558,237){\makebox(0,0)[lb]{\smash{{\mathmode{\scriptstyle 2}}}}}
\end{picture}
}}
\def\example
{{
\begin{picture}(7665,4095)(0,-10)
\thicklines
\path(825,3261)(2175,2061)(2175,261)
\path(2175,2061)(3525,3261)
\dottedline{45}(3975,486)(2400,1911)
\path(2509.112,1852.736)(2400.000,1911.000)(2468.857,1808.244)
\put(4050,336){\makebox(0,0)[lb]{\smash{{\mathmode{\scriptstyle f_{abc}\
(\text{or }f_{bca}\text{ or }f_{cab})}}}}}
\put(2250,1611){\makebox(0,0)[lb]{\smash{{\mathmode{\scriptstyle c}}}}}
\put(2325,2511){\makebox(0,0)[lb]{\smash{{\mathmode{\scriptstyle a}}}}}
\put(2250,561){\makebox(0,0)[lb]{\smash{{\mathmode{\scriptstyle c_1}}}}}
\put(1800,36){\makebox(0,0)[lb]{\smash{{\mathmode{\scriptstyle
f_{c_1c_2c_3}}}}}}
\put(3975,1311){\makebox(0,0)[lb]{\smash{{\mathmode{\scriptstyle
t^{c_3a_2}}}}}}
\put(1725,3936){\makebox(0,0)[lb]{\smash{{\mathmode{\scriptstyle t^{a_3b_2}\
(\text{or }t^{b_2a_3})}}}}}
\put(2925,2511){\makebox(0,0)[lb]{\smash{{\mathmode{\scriptstyle t^{aa_1}}}}}}
\put(1500,2736){\makebox(0,0)[lb]{\smash{{\mathmode{\scriptstyle t^{bb_1}}}}}}
\put(2175,2061){\ellipse{3600}{3600}}
\put(2325,1086){\makebox(0,0)[lb]{\smash{{\mathmode{\scriptstyle t^{cc_1}}}}}}
\put(1350,3486){\makebox(0,0)[lb]{\smash{{\mathmode{\scriptstyle b_2}}}}}
\put(675,1386){\makebox(0,0)[lb]{\smash{{\mathmode{\scriptstyle t^{b_3c_2}}}}}}
\put(1650,2136){\makebox(0,0)[lb]{\smash{{\mathmode{\scriptstyle b}}}}}
\put(825,2886){\makebox(0,0)[lb]{\smash{{\mathmode{\scriptstyle b_1}}}}}
\put(150,2811){\makebox(0,0)[lb]{\smash{{\mathmode{\scriptstyle b_3}}}}}
\put(1500,486){\makebox(0,0)[lb]{\smash{{\mathmode{\scriptstyle c_2}}}}}
\put(2775,186){\makebox(0,0)[lb]{\smash{{\mathmode{\scriptstyle c_3}}}}}
\put(2925,3111){\makebox(0,0)[lb]{\smash{{\mathmode{\scriptstyle a_1}}}}}
\put(3300,3636){\makebox(0,0)[lb]{\smash{{\mathmode{\scriptstyle a_3}}}}}
\put(3675,3336){\makebox(0,0)[lb]{\smash{{\mathmode{\scriptstyle f_{a_1a_2a_3}\
(\text{or }\ldots)}}}}}
\put(0,3411){\makebox(0,0)[lb]{\smash{{\mathmode{\scriptstyle
f_{b_1b_2b_3}}}}}}
\put(3900,2886){\makebox(0,0)[lb]{\smash{{\mathmode{\scriptstyle a_2}}}}}
\end{picture}
}}
\def\glNexample
{{
\begin{picture}(2509,2233)(0,-10)
\thicklines
\put(1137.500,972.500){\arc{1629.801}{1.4940}{3.9921}}
\put(1412.500,972.500){\arc{1629.801}{5.4327}{7.9308}}
\put(1401.275,952.092){\arc{1886.973}{5.4347}{7.9084}}
\put(1148.725,952.092){\arc{1886.973}{1.5164}{3.9901}}
\put(1275.000,1125.625){\arc{1818.750}{3.8759}{5.5488}}
\put(1275,985){\whiten\ellipse{150}{150}}
\put(1275,985){\ellipse{150}{150}}
\put(1275,85){\whiten\ellipse{150}{150}}
\put(1275,85){\ellipse{150}{150}}
\put(1950,1660){\whiten\ellipse{150}{150}}
\put(1950,1660){\ellipse{150}{150}}
\put(600,1660){\whiten\ellipse{150}{150}}
\put(600,1660){\ellipse{150}{150}}
\path(1200,910)(1200,160)
\path(1350,910)(1350,160)
\path(1275,1060)(1875,1660)
\path(1350,985)(1950,1585)
\put(1275.000,972.500){\arc{1825.000}{3.9948}{5.4299}}
\path(675,1660)(1275,1060)
\put(525,685){\makebox(0,0)[lb]{\smash{{\mathmode{\scriptstyle p}}}}}
\path(600,1585)(1200,985)
\put(1200,2110){\makebox(0,0)[lb]{\smash{{\mathmode{\scriptstyle u}}}}}
\put(1650,1135){\makebox(0,0)[lb]{\smash{{\mathmode{\scriptstyle i}}}}}
\put(1425,1435){\makebox(0,0)[lb]{\smash{{\mathmode{\scriptstyle j}}}}}
\put(975,1435){\makebox(0,0)[lb]{\smash{{\mathmode{\scriptstyle k}}}}}
\put(675,1135){\makebox(0,0)[lb]{\smash{{\mathmode{\scriptstyle l}}}}}
\put(1425,535){\makebox(0,0)[lb]{\smash{{\mathmode{\scriptstyle n}}}}}
\put(0,610){\makebox(0,0)[lb]{\smash{{\mathmode{\scriptstyle q}}}}}
\put(1950,685){\makebox(0,0)[lb]{\smash{{\mathmode{\scriptstyle r}}}}}
\put(2325,610){\makebox(0,0)[lb]{\smash{{\mathmode{\scriptstyle s}}}}}
\put(1200,1660){\makebox(0,0)[lb]{\smash{{\mathmode{\scriptstyle t}}}}}
\put(900,535){\makebox(0,0)[lb]{\smash{{\mathmode{\scriptstyle m}}}}}
\end{picture}
}}
\def\ijklmn
{{
\begin{picture}(1137,1221)(0,-10)
\thicklines
\path(525,537)(525,12)
\path(675,537)(675,12)
\path(600,687)(1050,1137)
\path(675,612)(1125,1062)
\path(150,1137)(600,687)
\put(600,612){\whiten\ellipse{150}{150}}
\put(600,612){\ellipse{150}{150}}
\path(75,1062)(525,612)
\put(225,162){\makebox(0,0)[lb]{\smash{{\mathmode{\scriptstyle m}}}}}
\put(975,762){\makebox(0,0)[lb]{\smash{{\mathmode{\scriptstyle i}}}}}
\put(750,1062){\makebox(0,0)[lb]{\smash{{\mathmode{\scriptstyle j}}}}}
\put(300,1062){\makebox(0,0)[lb]{\smash{{\mathmode{\scriptstyle k}}}}}
\put(0,762){\makebox(0,0)[lb]{\smash{{\mathmode{\scriptstyle l}}}}}
\put(750,162){\makebox(0,0)[lb]{\smash{{\mathmode{\scriptstyle n}}}}}
\end{picture}
}}
\def\soexample
{{
\begin{picture}(3616,3634)(0,-10)
\thicklines
\path(458,3012)(1808,1812)(1808,12)
\path(1808,1812)(3158,3012)
\put(1883,1062){\makebox(0,0)[lb]{\smash{{\mathmode{\scriptstyle c}}}}}
\put(1808,1812){\ellipse{3600}{3600}}
\put(2258,2562){\makebox(0,0)[lb]{\smash{{\mathmode{\scriptstyle a}}}}}
\put(1733,3312){\makebox(0,0)[lb]{\smash{{\mathmode{\scriptstyle f}}}}}
\put(908,2187){\makebox(0,0)[lb]{\smash{{\mathmode{\scriptstyle b}}}}}
\put(383,987){\makebox(0,0)[lb]{\smash{{\mathmode{\scriptstyle d}}}}}
\put(3158,987){\makebox(0,0)[lb]{\smash{{\mathmode{\scriptstyle e}}}}}
\end{picture}
}}
\begin{document}
\title{Lie Algebras and the Four Color Theorem}

\author{Dror Bar-Natan}
\address{Institute of Mathematics\\
	\AA rhus University\\
	8000 \AA rhus C \\
	Denmark}
\curraddr{Institute of Mathematics\\
        The Hebrew University\\
        Giv'at-Ram, Jerusalem 91904\\
        Israel}
\email{drorbn@math.huji.ac.il}

\thanks{This preprint is available electronically at
  {\tt http://www.ma.huji.ac.il/$\sim$drorbn}, at \newline
  {\tt file://ftp.ma.huji.ac.il/drorbn}, and at {\tt
    http://xxx.lanl.gov/abs/q-alg/9606016}.
}

\date{This edition: Jun.~23,~1996; \ \ First edition: August 18, 1995.}

\maketitle

\begin{abstract}
We present a ``reasonable'' statement about Lie algebras that is equivalent
to the Four Color Theorem.
\end{abstract}

\tableofcontents
\section{Introduction}

Let us start by recalling a well-known construction that associates
to any finite dimensional metrized Lie algebra $L$ a numerical-valued
functional $W_L$ defined on the set of all oriented trivalent graphs
$G$ (that is, trivalent graphs in which every vertex is endowed with
a cyclic ordering of the edges emanating from it).  This construction
underlies the gauge-group dependence of gauge theories in general
and of the Chern-Simons topological field theory in particular (see
e.g.~\cite{Bar-Natan:Thesis, AxelrodSinger:PCS-I, AxelrodSinger:PCS-II})
and plays a prominent role in the theory of finite type (Vassiliev)
invariants of knots (\cite{Bar-Natan:Weights, Bar-Natan:Vassiliev,
Bar-Natan:VasBib}) and most likely also in the theory of finite
type invariants of 3-manifolds (\cite{Ohtsuki:IntegralHomology,
GaroufalidisOhtsuki:3ManifoldsIII, Rozansky:RationalHomology}).

Fix a finite dimensional metrized Lie algebra $L$ (that is, a finite
dimensional Lie algebra with an {\it ad}-invariant symmetric non-degenerate
bilinear form),
choose some basis $\{L_a\}_{a=1}^{\dim L}$ of $L$, let $t_{ab}=\langle
L_a,L_b\rangle$ be the metric tensor, let $t^{ab}$ be the inverse matrix of
$t_{ab}$, and let
$f_{abc}$ be the structure constants of $L$ relative to $\{L_a\}$:
\[ f_{abc}=\langle L_a,[L_b,L_c]\rangle. \]
Let $G$ be some oriented trivalent graph. To define $W_L$, label all half-edges
of $G$ by symbols from the list $a,b,c,\ldots,a_1,b_1,\ldots$, and sum over
$a,b,\ldots,a_1,\ldots\in\{1,\ldots,\dim L\}$ the product over the vertices of
$G$ of the structure constants ``seen'' around each vertex times the
product over the edges of the $t^{\cdot\cdot}$'s seen on each edge. This
definition is much better explained by an example, as in
figure~\ref{example}.

\begin{figure}[htpb]
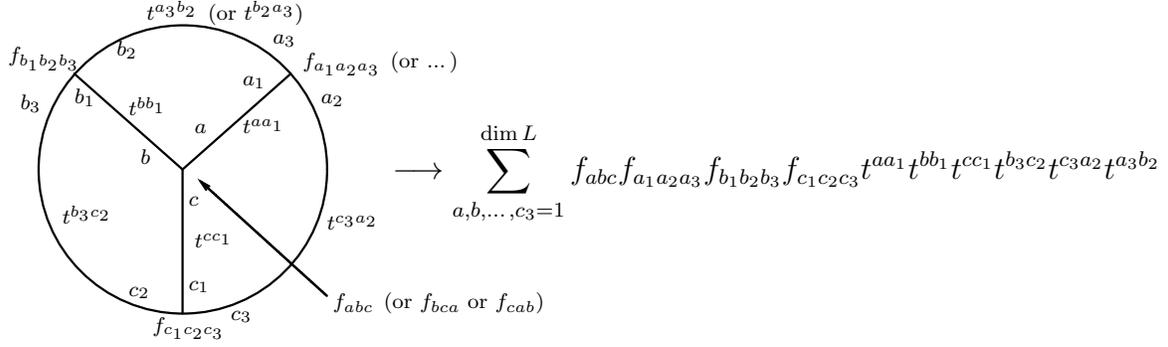

\[ \EEPIC{\example}{0.5}
  \hspace{-33mm} \longrightarrow
  \sum_{a,b,\ldots,c_3=1}^{\dim L}
    f_{abc}f_{a_1a_2a_3}f_{b_1b_2b_3}f_{c_1c_2c_3}
    t^{aa_1}t^{bb_1}t^{cc_1}t^{b_3c_2}t^{c_3a_2}t^{a_3b_2}
\]
\caption{An example illustrating the construction of $W_L(G)$. Notice that
when $G$ is drawn in the plane, we assume counterclockwise orientation for
all vertices (unless noted otherwise), and that the cyclic symmetry
$f_{abc}=f_{bca}=f_{cab}$ of the
structure constants and the symmetry $t^{ab}=t^{ba}$ of the inverse
metric ensures that $W_L(G)$ is well defined.}
\lbl{example} \end{figure}

By introducing an explicit change-of-basis matrix as
in~\cite{Bar-Natan:Weights} or by re-interpreting $W_L(G)$ in terms of
abstract tensor calculus as in~\cite{Bar-Natan:Vassiliev}, one can verify that
$W_L(G)$ does not depend on the choice of the basis $\{L_a\}$. Typically
one chooses a ``nice''  orthonormal (or almost orthonormal) basis $\{L_a\}$,
so that most of the constants $t^{ab}$ and $f_{abc}$ vanish, thus greatly
reducing the number of summands in the definition of $W_L(G)$.

Unless otherwise stated,
whenever dealing with a Lie algebra of matrices, we will
take the metric to be the matrix trace in the defining representation:
$\langle L_a,L_b\rangle = \tr(L_aL_b)$.

\begin{lemmadefinition} \lbl{LemmaDefinition} (proof in section~\ref{slN})
If a connected $G$ has $v$ vertices, then
$W_{sl(N)}(G)$ is a polynomial in $N$ of degree at most $\frac{v}{2}+2$ in
$N$. Thus we can set $W_{sl(N)}^{\text{top}}(G)$ to be the coefficient
of $N^{\frac{v}{2}+2}$ in $W_{sl(N)}(G)$.
\end{lemmadefinition}

The following statement sounds
rather reasonable; it just says that if $G$ is ``$sl(2)$-trivial'', then it
is at least ``$sl(N)$-degenerate''. For us who grew up thinking that all
that there is to learn about $sl(N)$ is already in $sl(2)$, this is not a
big surprise:

\begin{statement} \lbl{statement}
For a connected oriented trivalent graph $G$, $W_{sl(2)}(G)=0$ implies
$W_{sl(N)}^{\text{top}}(G)=0$.
\end{statement}

Lie-theoretically, there is much to say about $sl(2)$ and $sl(N)$. There
are representations of $sl(2)$ into $sl(N)$, there is an ``almost
decomposition'' of $sl(N)$ into a product of $sl(2)$'s\footnote{
See~\cite{Bar-NatanGaroufalidis:MMR} for a similar context in which the
different $sl(2)$'s ``decouple''.}, and there are many
other similarities. A-priori, the above statement sounds within reach. The
purpose of this note is to explain why statement~\ref{statement} is
equivalent to the Four Color Theorem\footnote{The Four Color Theorem
was conjectured by
Francis Guthrie in 1852 and proven by K.~I.~Appel and
W.~Haken~\cite{AppelHaken:Book} in 1976. See
also~\cite{SaatyKainen:4CP}.}.

This equivalence follows from the following two propositions, proven in
sections~\ref{slN} and~\ref{sl2}, respectively:

\begin{proposition} \lbl{slNprop} Let $G$ be a connected oriented trivalent
graph. If $G$ is 2-connected,
$|W_{sl(N)}^{\text{top}}(G)|$ is is equal to the number
of embeddings of $G$ in an oriented sphere. Otherwise,
$W_{sl(N)}^{\text{top}}(G)=0$.
\end{proposition}

\begin{proposition} \lbl{sl2prop} (Penrose~\cite{Penrose:Tensors}. See
also~\cite{Kauffman:CrossProduct, Kauffman:SpinNets,
KauffmanSaleur:PlanarColoring}.)
If $G$ is planar with $v$ vertices and $G^c$ is
the map defined by its
complement, than $|W_{sl(2)}(G)|$ is $2^{\frac{v}{2}-2}$ times
the total number of colorings of $G^c$ with four colors so that adjacent
states are colored with different colors.
\end{proposition}

Indeed, statement~\ref{statement} is clearly equivalent to
\[ |W_{sl(N)}^{\text{top}}(G)|\neq 0
  \quad\Rightarrow\quad
  |W_{sl(2)}(G)|\neq 0,
\]
which by propositions~\ref{slNprop} and~\ref{sl2prop} is the same as saying
\[ \left(\parbox{1.9in}{
    \begin{center}
      $G$ has a planar embedding with $G^c$ a map
    \end{center}
  }\right)
  \quad\Rightarrow\quad
  \left(\parbox{1.4in}{$G^c$ has a 4-coloring}\right).
\]
Notice that if $G$ is connected, $G^c$ is a map (does not have states that
border themselves) iff $G$ is 2-connected.

\begin{remark}
We've chosen the formulation of statement~\ref{statement} that we felt was
the most appealing. With no change to the end result, one can replace
$sl(N)=A_{N-1}$ by $B_N$, $C_N$, $D_N$, or $gl(N)$ and $sl(2)$ by $so(3)$
in the formulation of statement~\ref{statement}. In fact, in the proofs we
actually work with $gl(N)$ and $so(3)$ rather than with $sl(N)$ and
$sl(2)$.
\end{remark}

\subsection{Acknowledgement}
I wish to thank D.\ Goldberg, J.\ Goldman,
M.\ Hutchings, L.\ Kauffman, D.\ P.\ Thurston
and M.\ Wunderlich for many helpful conversations.
\section{Understanding $W_{sl(N)}$} \lbl{slN}

As Lie algebras, $gl(N)$ is just $sl(N)$ plus an Abelian factor. As Abelian
Lie algebras have vanishing structure constants,
$W_{sl(N)}(G)=W_{gl(N)}(G)$ for any oriented trivalent graph $G$. So let
us concentrate on computing $W_{gl(N)}(G)$ for such $G$. For the basis of
$gl(N)$, we pick the matrices $\{L_a\}_{a=1}^{N^2}=\{L_{ij}\}_{i,j=1}^N$,
where $L_{ij}$ is the matrix with $1$ in the $ij$ entry and $0$ everywhere
else. As the basis is indexed by a double index rather than by a single
index, it is convenient to label every half-edge of $G$ by two symbols from
the list $i,j,\ldots,i_1,\ldots$ and double all the edges:
\begin{equation} \lbl{DoubleEdges} \EEPIC{\DoubleEdges}{0.5}. \end{equation}
The metric $t_{ab}=t_{(ij)(kl)}$ of $gl(N)$ is given by $t_{(ij)(kl)}=\tr
L_{ij}L_{kl}=\delta_{jk}\delta_{il}$, and its inverse is given by the same
formula:
\[ t^{(ij)(kl)}=\delta_{jk}\delta_{li}. \]
This formula means that in the summation defining $W_{gl(N)}(G)$ we can
assume the equalities $j=k$ and $l=i$ along each edge as
in~\eqref{DoubleEdges}. In other words, it is enough to label every
doubled edge with just one pair of indices, getting an overall picture like
\[ \EEPIC{\glNexample}{0.5} \longrightarrow
  \sum_{i,j,\ldots,u=1}^N
  f_{(ij)(kl)(mn)}f_{(ji)(rs)(ut)}f_{(lk)(tu)(qp)}f_{(nm)(pq)(sr)},
\]
where $f_{(ij)(kl)(mn)}$ are the structure constants in our basis:
\[ \EEPIC{\ijklmn}{0.5} =
  f_{(ij)(kl)(mn)} =
  \langle L_{ij}, [L_{kl},L_{mn}] \rangle =
  \tr(L_{ij}L_{kl}L_{mn})-\tr(L_{mn}L_{kl}L_{ij})
\]
\begin{equation} \lbl{ThickVertex}
  = \delta_{jk}\delta_{lm}\delta_{ni}-\delta_{nk}\delta_{li}\delta_{jm} =
  \EEPIC{\ThickVertex}{0.5}.
\end{equation}
In the last equation, indices connected by a line can be assumed to be
equal in the summation defining $W_{gl(N)}(G)$. Once the edges and vertices
of $G$ are ``thickened'' as in~\eqref{DoubleEdges}
and~\eqref{ThickVertex}, the summation over $i,j,\ldots$ becomes the
counting of the number of solutions of the equalities determined by the
connected components of the thickened~$G$. This number is simply $N$ raised
to the number of connected components:
\begin{multline*}
  \EEPIC{\ThickG}{0.5} \longrightarrow
  \sum_{i,j,\ldots,u=1}^N \delta_{il}\delta_{lp}\cdots
  \\
  =\#\left\{
    1\leq i,j,\ldots\leq N:
    \begin{array}{c}
      i=l=p=m=j=t=k=n=r \\
      u=q=s
    \end{array}
  \right\}
  =N^2.
\end{multline*}

Summarizing, we find the formula\footnote{Compare with \cite[equation
(36)]{Bar-Natan:Vassiliev}; for similar formulas in the cases of $so(N)$
and $sp(N)$, see \cite[equation (33)]{Bar-Natan:Vassiliev} and
\cite[exercise 6.37]{Bar-Natan:Vassiliev}.}
\begin{equation} \lbl{WglN}
  W_{gl(N)}(G)=\sum_{\text{markings $M$ of $G$}} \sign(M) N^{b(T_M)},
\end{equation}
where:
\begin{itemize}
\item A markings $M$ of $G$ is a marking of each vertex of $G$ by a sign in
  $\{+,-\}$, and  $\sign(M)$ is the product of these signs.
\item The thickening $T_M$ corresponding to a marking $M$ is the oriented
  surface with boundary obtained from $G$ as follows:
  \begin{itemize}
  \item Replace the vertices marked by a ``$+$'' with ``joints'' and the
    vertices marked by a ``$-$'' with ``twisted joints'' as
    in~\eqref{ThickVertex}.
  \item Orient these surface pieces using
    the following black($b$)/white($w$)
    convention for the thickening of vertices:
     \[ \EEPIC{\BWConvention}{0.75}. \]
    In other words, ``+''-vertices are embedded in the thickening of $G$
    so that they are seen as oriented counterclockwise from the white
    side of the thickening, while ``-''-vertices are seen as oriented
    clockwise from the white side of the thickening.
  \item Finally, connect the joints together along the edges of $G$ by
    bands, in the only way consistent with the orientations of the joints.
  \end{itemize}
\item $b(T_M)$ is the number of boundary components of $T_M$.
\end{itemize}

If $M$ is a marking of $G$, let $S_M$ be the closed oriented
surface obtained by
gluing a disk into each boundary component of the thickening $T_M$.
With $\chi$ denoting Euler characteristic and $g$ denoting genus, we have
\[ 2-2g(S_M)=\chi(S_M)=\chi(T_M)+b(T_M)=\chi(G)+b(T_M). \]
Remembering that $G$ is trivalent and thus $\chi(G)=-\frac{v}{2}$, we get
\[ b(T_M)=-\chi(G)+2-2g(S_M)=\frac{v}{2}+2-2g(S_M). \]
Thus $b(T_M)$ is maximal when $g(S_M)=0$ and in that case
$b(T_M)=\frac{v}{2}+2$. With~\eqref{WglN}, this proves
lemma~\ref{LemmaDefinition}. Furthermore,
calling a marking $M$ {\em spherical} when $S_M$ is a sphere, we get the
formula
\[ W_{sl(N)}^{\text{top}}(G) = W_{gl(N)}^{\text{top}}(G) =
  \sum_{\text{spherical markings $M$ of $G$}} \sign(M).
\]

\begin{pf*}{Proof of proposition~\ref{slNprop}} Let $G$ be 2-connected, and
consider the map
\[ \Theta:\{\text{spherical markings $M$ of $G$}\}
   \longrightarrow
   \{\text{embeddings of $G$ in an oriented sphere}\}
\]
defined by mapping $M$ to the natural embedding of $G$ in $S_M$.
It is clear that $\Theta$ is a bijection. Indeed, if an
embedding of $G$ in an oriented sphere $S^2$ is given, one can reconstruct
$M$ by
marking the vertices that are oriented counterclockwise within $S_M$
(as seen from the white side of $S_M$) by a ``$+$'' sign, and marking
all other vertices by a ``$-$''.

To conclude the proof, it is enough to show that for spherical markings,
$\sign(M)$ is independent of the marking. Clearly,
$\sign(M)$ is also equal to $(-1)^{v_-}$, where $v_-$ is the
number of vertices of $G$ that are embedded clockwise in $S_M$ by
$\Theta(M)$ (with $S_M$ viewed from its white side). By a theorem of
H.~Whitney~\cite{Whitney:TopologicalInvariants,
Whitney:Classification}\footnote{Check~\cite[corollary
6]{Lipson:LinkSignature} for a version closer to what we need, and remember
that our graph is 2-connected and hence some of the moves
in~\cite{Lipson:LinkSignature} are irrelevant for us.},
one can get from any spherical embedding of
a trivalent 2-connected graph $G$ to any other such embedding by a
sequence of
flips as in figure~\ref{Flip}. Such flips do not change the parity of
$v_-$, since the number of vertices that is flipped is even.
\end{pf*}

\begin{figure}[htpb]
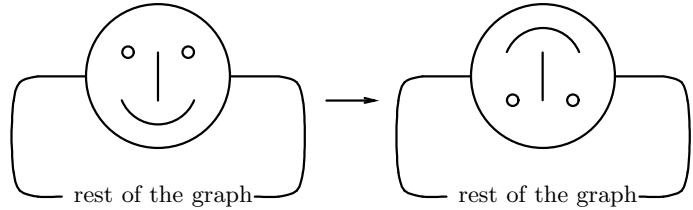

\[ \EEPIC{\Flip}{0.5} \]
\caption{A flip takes a part of a graph that connects to the rest via only
two edges, and flips it over.} \lbl{Flip}
\end{figure}
\section{Understanding $W_{sl(2)}$} \lbl{sl2}

Proposition~\ref{sl2prop} is due to Penrose~\cite{Penrose:Tensors} (see
also~\cite{Kauffman:CrossProduct, Kauffman:SpinNets,
KauffmanSaleur:PlanarColoring}). For completeness, we reproduce its proof
in this section.

As Lie algebras, $sl(2)$ is isomorphic to $so(3)$, so let us work with
$so(3)$ instead. The standard basis of $so(3)$ is given by the matrices
\[
  L_1 = \left(
    \begin{array}{rrr} 0 & 0 & 0 \\ 0 & 0 & -1 \\ 0 & 1 & 0 \end{array}
  \right), \quad
  L_2 = \left(
    \begin{array}{rrr} 0 & 0 & 1 \\ 0 & 0 & 0 \\ -1 & 0 & 0 \end{array}
  \right)
  \quad \text{and} \quad
  L_3 = \left(
    \begin{array}{rrr} 0 & -1 & 0 \\ 1 & 0 & 0 \\ 0 & 0 & 0 \end{array}
  \right).
\]
Let us pick the scalar product $\langle\cdot,\cdot\rangle$
on $so(3)$ to be the one that makes this
basis orthonormal, and let us denote the corresponding functional on graphs
by $\Wst$, with the ``$\widetilde{\quad}$'' on top of the $so(3)$ to remind
us that we are not using the standard matrix-trace scalar product.

One can easily verify that $\frac{1}{2}\langle\cdot,\cdot\rangle$ is the
scalar product induced on $so(3)$ from matrix-trace in $sl(2)$. Thus,
remembering that in the construction of $W_L$ vertices scale with the
scalar product and edges scale with its inverse, we find that
\begin{equation} \lbl{sl2so3}
  W_{sl(2)}(G) = \left(\frac{1}{2}\right)^{v-e} \Wst(G)
  =2^{\frac{v}{2}} \Wst(G),
\end{equation}
where $G$ is an oriented trivalent graph with $v$ vertices and $e$ edges.
Let us fix such a $G$ once and for all, and let us assume that it is planar
and that all the vertices of $G$ are oriented counterclockwise in the
plane. Flipping the orientation of any given vertex just reverses the sign
of $\Wst(G)$, and so the latter assumption does not limit the generality of
our arguments.

\begin{pf*}{Proof of proposition~\ref{sl2prop}} With~\eqref{sl2so3} in
mind, proposition~\ref{sl2prop} clearly follows from the following two
lemmas.
\end{pf*}

\begin{lemma} \lbl{PenroseLemma} (Penrose~\cite{Penrose:Tensors}. See
also~\cite{Kauffman:CrossProduct, Kauffman:SpinNets,
KauffmanSaleur:PlanarColoring}.)
For a planar $G$ as above,
$|\Wst(G)|$ is the number of colorings of the edges of $G$
with three colors $\{1,2,3\}$, so that the edges emanating from any single
vertex are of different colors.
\end{lemma}

\begin{lemma} \lbl{TaitTheorem} (Tait's theorem~\cite{Tait:Colouring})
Edge-3-colorings as in the previous lemma are in a bijective
correspondence with 4-colorings of the map
$G^c$ that fix the color of the ``state at infinity''.
\end{lemma}

\begin{pf*}{Proof of lemma~\ref{PenroseLemma}}
In the basis $\{L_a\}$, the structure constants of $so(3)$ are given by
\[ f_{abc}=\epsilon_{abc}=\begin{cases}
    \sign(abc)	& \text{if $abc$ is a permutation,} \\
    0		& \text{otherwise.}
  \end{cases}
\]
Remembering also that $\{L_a\}$ is orthonormal by choice, the computation
of $\Wst(G)$ is given (on a simple example) by:
\[ \EEPIC{\soexample}{0.4} \longrightarrow
  \sum_{a,b,c,d,e,f=1}^3
  \epsilon_{abc}\epsilon_{aef}\epsilon_{bfd}\epsilon_{cde}.
\]
The $\epsilon$ symbols force the indices coming into each vertex to be
different, and hence clearly
\begin{equation} \lbl{SPV}
  \Wst(G) = \sum_{\text{edge-3-colorings of }G}
  \prod\left(\text{a sign per vertex}\right),
\end{equation}
where the sign at
each vertex is the sign of the permutation of $\{1,2,3\}$ induced by an
edge-3-coloring, as read counterclockwise around the vertex.

The only thing left to show is that the product of signs in~\eqref{SPV} is
independent of the edge-3-coloring. A clever way to do that, discovered by
L.~H.~Kauffman, is to replace every edge colored by a ``3'' by a pair of
edges colored ``1'' and ``2'' (in symbols, $1+2=3$). This defines two
families of circles in the plane, labeled by ``1'' and by ``2'' (see
figure~\ref{TwoFamilies}). By lumping together the signs on each end of a
``3'' edge and taking the product over all of those edges, one sees that
the overall sign depends only on the parity of the number of ``3'' edges
(always $\frac{e}{3}$), and the $\bold Z/2\bold Z$ intersection number of
the ``1'' family of circles with the ``2'' family of circles. By the Jordan
curve theorem, the latter is always $0$.
\end{pf*}

\begin{figure}[htpb]
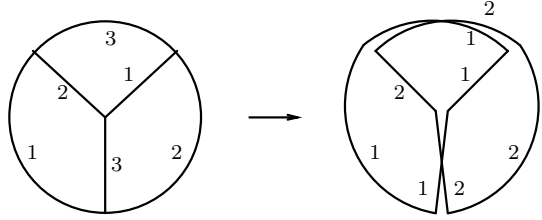

\[ \EEPIC{\TwoFamilies}{0.5} \]
\caption{The two families of circles obtained by splitting every ``3''
  edge.}
\lbl{TwoFamilies} \end{figure}

\begin{pf*}{Proof of lemma~\ref{TaitTheorem}}
This is a well known result (see
e.g.~\cite[theorem 9.12]{BondyMurty:GraphTheory}),
so let us only sketch the proof. Consider the
group $H=\bold Z/2\bold Z\times\bold Z/2\bold Z$. Given any 4-coloring of
$G^c$ by elements of $H$, one may associate to it an
edge-3-coloring
of $G$ by the non-zero elements of $H$, by coloring every edge by
the difference of the colors in the two faces adjacent to it. One then
verifies that this edge 3-coloring is well defined and that we get
a bijection
between the set of 4-colorings of $G^c$ that color the state at infinity
with $0$ and the set of edge-3-coloring of $G$.
\end{pf*}
\ifx\undefined\bysame
        \newcommand{\bysame}{\leavevmode\hbox to3em{\hrulefill}\,}
\fi

\end{document}